\documentclass[12pt]{amsart}

\usepackage{hyperref}
\usepackage{amsmath, amsthm, amssymb}
\usepackage{url} 
\usepackage{graphicx}
\usepackage{xcolor}
\usepackage{moreverb}

\usepackage{fancyvrb}

\newcommand{\R}{\mathbb R}
\newcommand{\V}{\mathcal V}
\newcommand{\X}{\mathfrak X}

\usepackage{listings}
 \lstnewenvironment{code}[1][]
  {\lstset{#1}}
  {}

\newtheorem{theorem}{Theorem}[section]
 \newtheorem{proposition}[theorem]{Proposition}
 \newtheorem{corollary}[theorem]{Corollary}
\theoremstyle{definition}
\newtheorem{definition}[theorem]{Definition}
\newtheorem{example}[theorem]{Example}
\newtheorem{remark}[theorem]{Remark}
\newtheorem{lemma}[theorem]{Lemma}

\begin{document}

\title[Surfaces of constant sectional curvature]{Surfaces of constant sectional curvature in Minkowski space with a canonical semi-symmetric connection}

\author{Muhittin Evren Aydin}
\address{Department of Mathematics, Faculty of Science, Firat University, 23200 Elazig,   Turkey}
\email{meaydin@firat.edu.tr}

\author{Rafael L\'opez}
\address{ Department of Geometry and Topology\\  University of Granada, 18071 Granada, Spain}
\email{rcamino@ugr.es}

\begin{abstract}
We study surfaces in Minkowski space $\mathbb{R}^3_1$ with constant intrinsic sectional curvature with respect
to canonical semi-symmetric connections, metric or non-metric, determined by a parallel vector field $\mathcal{V}$. We compute the ambient sectional curvature and obtain a formula for the intrinsic sectional curvature of nondegenerate surfaces in terms of their Gauss and mean curvatures. We provide a complete classification of cylindrical surfaces with constant intrinsic sectional curvature, covering all causal characters of $\mathcal{V}$ and of the rulings. For rotational surfaces whose axes are parallel to $\mathcal{V}$, we reduce the constant-curvature condition to an ordinary differential equation. Explicit solutions are obtained for null axes and, in the metric case, for surfaces with zero intrinsic sectional curvature.
\end{abstract}
\maketitle

\section{Introduction and formulation of the problem} \label{intro}

Let $(M,\langle\,,\,\rangle)$ be a pseudo-Riemannian manifold and let
$\nabla^0$ denote its Levi-Civita connection. Besides the Levi-Civita
connection, several affine connections have been introduced
to investigate geometric properties. Semi-symmetric connections, introduced in \cite{fs}, form an
important class. See also \cite{ds,ya}.

The torsion tensor $T^0$ of $\nabla^0$ vanishes identically. An affine connection $\nabla$ on $M$
is called a semi-symmetric connection \cite{fs} if its torsion tensor $T$, defined by $T(X,Y)=\nabla_{X}Y-\nabla_{Y}X-[X,Y]$, satisfies
\begin{equation*}
T(X,Y)=\langle\V,Y\rangle X-\langle\V,X\rangle Y,
\end{equation*}
for all $X,Y\in\X(M)$, where $\V$ is a nonzero vector field on $M$. In addition, $\nabla$ is called a metric connection if the metric $\langle\,,\,\rangle$ is parallel with respect to $\nabla$ \cite{hay}. Otherwise, $\nabla$ is called a non-metric connection. Semi-symmetric non-metric connections were introduced by Agashe   in \cite{ag}. See also Agashe and Chafle\cite{ac}.

Geometric properties of submanifolds of (pseudo-)Riemannian manifolds endowed with semi-symmetric metric and non-metric connections have been investigated by several authors. See, for example, \cite{im,lyl,mc1,mo,mmme,na,zz}, and \cite{mma,yay0,yay1} in the pseudo-Riemannian setting.

To unify the metric and non-metric cases, let $\alpha\in\{0,1\}$ and
consider the semi-symmetric connection
\begin{equation}\label{con0}
\nabla_XY=\nabla^0_XY+\langle \V, Y\rangle X- \alpha \langle X, Y\rangle \V, \quad X,Y\in\X(M),
\end{equation}
where $\V\in\X(M)$ is a nonzero vector field. This connection is metric for $\alpha=1$ and non-metric for $\alpha=0$.

The connection \eqref{con0} is a special case of the connection given in \cite[Equation (2.4)]{tr}, obtained by setting
$$\varphi_1=-\mathrm{Id}, \quad
\varphi_2=0, \quad
U=-U_2=-\V, \quad
f_1=0, \qquad
f_2=\alpha-1.
$$

Let $R$ denote the curvature tensor of the semi-symmetric connection
$\nabla$ on $M$ defined by \eqref{con0}. We write $R(X,Y;Z,W)=\langle R(X,Y)Z,W\rangle$. Let $\Pi\subset T_pM$ be a nondegenerate $2$-plane at
$p\in M$ and $\{e_1,e_2\}$ a pseudo-orthonormal basis of $\Pi$, with $\varepsilon_i=\langle e_i,e_i\rangle\in\{-1,1\}$. If
$\nabla$ is a semi-symmetric metric connection, then the quantity
\begin{equation}\label{sec0}
R(e_1,e_2;e_2,e_1)
:=
\langle R(e_1,e_2)e_2,e_1\rangle
\end{equation}
defines, after dividing by $\varepsilon_1\varepsilon_2$, the sectional curvature of the plane $\Pi$. However, if
$\nabla$ is non-metric, then the quantity
\eqref{sec0} does not define a sectional curvature, because it depends on
the choice of the orthonormal basis of $\Pi$.

To overcome this difficulty, Mihai and Mihai \cite{mm0} introduced the following notion of sectional curvature for semi-symmetric non-metric connections in the Riemannian setting. This notion was also used in the Lorentzian setting; see \cite{mma}.

\begin{definition}[\cite{mm0,mma}]\label{defK}
Let $\nabla$ be the connection defined by \eqref{con0}. Let $\Pi\subset T_pM$ be a nondegenerate plane and $\{e_1,e_2\}$ a pseudo-orthonormal basis of $\Pi$, with $\varepsilon_i=\langle e_i,e_i\rangle$. The sectional curvature of $\Pi$ with respect to $\nabla$ is defined by
\begin{equation}\label{def1}
K(\Pi)
=
\frac{1}{2\varepsilon_1\varepsilon_2}
\left[
R(e_1,e_2;e_2,e_1)
+
R(e_2,e_1;e_1,e_2)
\right].
\end{equation}
\end{definition}

In the Riemannian case ($\varepsilon_1=\varepsilon_2=1$) this definition coincides with that of \cite{mm0}, where it is proved that $K(\Pi)$ does not depend on the orthonormal basis of $\Pi$. In the Lorentzian case, the orthonormal bases of a timelike plane are related by Lorentz boosts, a situation not covered by \cite{mm0}.  However, for the connections considered in this paper, the definition remains independent of the choice of pseudo-orthonormal basis. This follows from the explicit expressions of Propositions~\ref{pr23} and \ref{pr31}, which only involve quantities that do not depend on the basis.

Using the notion of sectional curvature defined by \eqref{def1}, for $\alpha=0$, and using the notion of sectional curvature defined by  \eqref{def1}, the
present authors, together with Mihai, classified the surfaces of Euclidean space $\R^3$ with 
constant sectional curvature that are
invariant under one-parameter groups of translations and rotations
\cite{alm0,alm1}.

Notice that the definition \eqref{def1} is valid for both values of $\alpha\in\{0,1\}$. When $\alpha=1$, the connection $\nabla$ is
metric and therefore its curvature tensor satisfies the skew-symmetry property $R(X,Y;Z,W)=-R(X,Y;W,Z)$. Combining this with $R(Y,X)=-R(X,Y)$, we derive $R(e_2,e_1;e_1,e_2)=-R(e_1,e_2;e_1,e_2)=R(e_1,e_2;e_2,e_1)$. Consequently,
$$
K(\Pi)
=
\frac{1}{2\varepsilon_1\varepsilon_2}
\left[
R(e_1,e_2;e_2,e_1)
+
R(e_2,e_1;e_1,e_2)
\right]
=
\frac{R(e_1,e_2;e_2,e_1)}{\varepsilon_1\varepsilon_2},
$$
which coincides with the sectional curvature in the case $\alpha=1$.

For each fixed $\alpha\in\{0,1\}$, the semi-symmetric connection $\nabla$ given by \eqref{con0} is determined by the choice of the vector field $\V$. If $\V=0$, then $\nabla$ reduces to the Levi-Civita connection $\nabla^0$. The simplest nontrivial case is when $\V$ is  parallel with respect to $\nabla^0$. This motivates the following definition.

\begin{definition}\label{cano}
A semi-symmetric connection $\nabla$ defined by \eqref{con0} is said to be canonical if its defining vector field $\V$ is $\nabla^0$-parallel, that is, $\nabla^0\V=0$. Such a connection is called a canonical semi-symmetric connection (css-connection) determined by $\V$. In particular, if $\nabla$ is metric (resp. non-metric), then it is called a metric (resp. non-metric) css-connection determined by $\V$.
\end{definition}

Let the ambient space be the Minkowski $3$-space $\R^3_1$, that is, the affine space
$\R^3$ endowed with the Lorentzian metric $\langle\,,\,\rangle=dx^2+dy^2-dz^2,$ where $(x,y,z)$ are the Cartesian coordinates. Following the ideas in \cite{alm0,alm1}, the aim of the present paper is to classify the surfaces of $\R^3_1$ that are 
 invariant under a one-parameter group of isometries (namely, cylindrical   and rotational surfaces) and have constant sectional curvature with respect to a css-connection. Translation surfaces of $\R^3_1$ that are minimal with respect to semi-symmetric connections were studied in \cite{wa}. Here, the condition is of a different nature, as it concerns the sectional curvature \eqref{def1}. 

The main results are the following. In Section~\ref{sec2}, we compute the sectional curvature of $\R^3_1$ with respect to a css-connection (Proposition~\ref{pr23}). In Section~\ref{sec3}, we obtain a formula for the sectional curvature $K^\Sigma$ of a nondegenerate surface $\Sigma$ in terms of its Gauss and mean curvatures (Proposition~\ref{pr31}). This formula can also be expressed intrinsically as
$$
K^\Sigma=G-\frac{1+\alpha}{2}\,\Delta^\Sigma \xi+\frac{(1-\alpha)^2}{2}\,\langle\nabla^{0_\Sigma}\xi,\nabla^{0_\Sigma}\xi\rangle,\qquad \xi=\langle\V,\cdot\rangle,
$$
as shown in Corollary~\ref{intrinsic}. In Section~\ref{sec4}, we show that, for a cylindrical surface with generating curve $\gamma$, $K^\Sigma$ depends only on $p=\langle\V,\gamma'\rangle$ and $p'$. We classify all cylindrical surfaces with constant $K^\Sigma$, considering every causal character of $\V$ and of the rulings (Theorem~\ref{clascyl} and Proposition~\ref{nullrulings}). In Section~\ref{sec5}, we study rotational surfaces whose axis is parallel to $\V$. We prove that the condition $K^\Sigma=K_0$ reduces to an ordinary differential equation for the profile curve (Proposition~\ref{throt}). We explicitly integrate this equation when the axis is null (Theorem~\ref{thnull}). For the metric connection, we also obtain an explicit integration in the case of flat surfaces (Corollary~\ref{rotflat}).

\section{Curvature of Minkowski $3$-space with respect to css-connections} \label{sec2}

Let $(M,\langle\,,\,\rangle)$ be a pseudo-Riemannian manifold and $\nabla^0$ its Levi-Civita connection. Let $\V$ be  a nonzero vector field on $M$ and set $\delta:=\langle\V,\V\rangle\in\{-1,0,1\}$. Define
$$
\eta(X):=\langle\V,X\rangle,
$$
where $\V$ is assumed to be normalized when it is non-null (see Remark~\ref{rem-norm}). We also define the $(1,2)$-tensor field $B$ by 
\begin{equation}\label{Btens}
B=I\otimes \eta-\alpha\langle \, , \, \rangle\otimes\V,
\end{equation}
where $I$ denotes the identity tensor field on $TM$ and the tensor products are understood so that $B(X,Y)=\eta(Y)X-\alpha \langle X,Y\rangle\V.$

If $\V$ is $\nabla^0$-parallel, then $B$ is also $\nabla^0$-parallel. Indeed, since $\nabla^0$ is the Levi-Civita connection, we have $\nabla^0\langle \, , \, \rangle=0$ and $\nabla^0I=0$. Moreover, $\nabla^0\V=0$ implies $\nabla^0\eta=0$. Thus, by the Leibniz rule,
$$
\nabla^0B
=(\nabla^0 I)\otimes\eta
+I\otimes(\nabla^0\eta)-\alpha(\nabla^0\langle \, , \, \rangle)\otimes\V
-\alpha\langle \, , \, \rangle\otimes(\nabla^0\V)
=0.
$$

We first express the curvature tensor of a css-connection in terms of the curvature tensor of the Levi-Civita connection and its defining vector field $\V$.

\begin{lemma}\label{lem2}
The curvature tensor $R$ of a css-connection determined by $\V$
satisfies
\begin{equation}\label{curvature}
\begin{aligned}
R(X,Y)Z
={}&R^0(X,Y)Z+\bigl(\eta(Y)\eta(Z)-\alpha\delta
\langle Y,Z\rangle\bigr)X\\
&-\bigl(\eta(X)\eta(Z)-\alpha\delta
\langle X,Z\rangle\bigr)Y\\
&+\alpha^2
\bigl(
\eta(X)\langle Y,Z\rangle
-\eta(Y)\langle X,Z\rangle
\bigr)\mathcal V,
\end{aligned}
\end{equation}
where $R^0$ is the curvature tensor of $\nabla^0$ and $\delta=\langle \V,\V\rangle$.
\end{lemma}
\begin{proof}
The css-connection given by \eqref{con0} can be written as
\[
\nabla_XY=\nabla^0_XY+B(X,Y),
\qquad X,Y\in\X(M),
\]
where $B$ is the $(1,2)$-tensor field defined by \eqref{Btens}. A standard computation yields
$$
\begin{aligned}
R(X,Y)Z
={}&R^0(X,Y)Z
+(\nabla_X^0B)(Y,Z)
-(\nabla_Y^0B)(X,Z)\\
&+B(X,B(Y,Z))-B(Y,B(X,Z)).
\end{aligned}
$$
Since $B$ is $\nabla^0$-parallel, we have
\begin{equation}\label{curvtensors}
R(X,Y)Z
=R^0(X,Y)Z+B(X,B(Y,Z))-B(Y,B(X,Z)).
\end{equation}
Next,
\begin{align*}
B(X,B(Y,Z))
={}&\eta(B(Y,Z))X
-\alpha\langle X,B(Y,Z)\rangle\mathcal V.
\end{align*}
Using $\eta(\mathcal V)=\delta$, we obtain
\[
\eta(B(Y,Z))
=
\eta(Y)\eta(Z)-\alpha\delta\langle Y,Z\rangle
\]
and
\[
\langle X,B(Y,Z)\rangle
=
\eta(Z)\langle X,Y\rangle
-\alpha\eta(X)\langle Y,Z\rangle.
\]
Consequently,
\begin{align*}
B(X,B(Y,Z))
={}&
\bigl(\eta(Y)\eta(Z)-\alpha\delta\langle Y,Z\rangle\bigr)X\\
&-\alpha\eta(Z)\langle X,Y\rangle\V
+\alpha^2\eta(X)\langle Y,Z\rangle\V.
\end{align*}
Similarly,
\begin{align*}
B(Y,B(X,Z))
={}&
\bigl(\eta(X)\eta(Z)-\alpha\delta\langle X,Z\rangle\bigr)Y\\
&-\alpha\eta(Z)\langle X,Y\rangle\V
+\alpha^2\eta(Y)\langle X,Z\rangle\V.
\end{align*}
Subtracting these identities, the terms
$-\alpha\eta(Z)\langle X,Y\rangle\V$ cancel.
Substituting the resulting expression into
\eqref{curvtensors} proves \eqref{curvature}.
\end{proof}

From now on, let the ambient space $M=\R^3_1$ be the Minkowski $3$-space endowed with the indefinite metric $\langle \, , \, \rangle=dx^2+dy^2-dz^2$. Since every $\nabla^0$-parallel vector field on $\R^3_1$ is constant, the defining vector field
$\V$ is constant. Moreover, since the curvature tensor $R^0$ of $\R^3_1$ vanishes identically, formula \eqref{curvature} reduces to
\begin{equation}\label{ss}
\begin{split}
R(X,Y)Z&=\bigl(\eta(Y)\eta(Z)-\alpha\delta\langle Y,Z\rangle\bigr)X
-\bigl(\eta(X)\eta(Z)-\alpha\delta\langle X,Z\rangle\bigr)Y
\\
&+\alpha^2\bigl(\eta(X)\langle Y,Z\rangle-\eta(Y)\langle X,Z\rangle\bigr)\V.
\end{split}
\end{equation}

\begin{remark}[Normalization of $\V$]\label{rem-norm}
Let $\Phi$ be a linear isometry of $\R^3_1$. Then $\Phi$ carries the css-connection determined by
$\V$ to the css-connection determined by $\Phi(\V)$,
preserving the corresponding sectional curvatures. For $\lambda>0$, let $\Phi_\lambda(q)=\lambda q$ be a homothety. If $\Sigma$ is a nondegenerate surface, then, under $\Phi_\lambda$, its Gauss   and mean curvature transform as $G\mapsto\lambda^{-2}G$ and
$H^0\mapsto\lambda^{-1}H^0$, respectively, while the unit normal $N$ remains unchanged at corresponding points. Hence, by \eqref{sectcurvat} below, the sectional curvature of $\Phi_\lambda(\Sigma)$ with respect to the css-connection determined by $\V/\lambda$
is $\lambda^{-2}$ times that of $\Sigma$ with respect to the css-connection determined by $\V$, at corresponding points. Consequently, the classification of surfaces with constant sectional curvature for $\V$ is equivalent, up to homotheties and rescaling of the curvature constant, to that for any positive multiple of $\V$. Negative multiples are handled by the isometry $-\mathrm{Id}$. This justifies the normalization of $\V$
whenever it is non-null. When $\V$ is null, for every $\lambda>0$, there exists a Lorentz boost $\mathcal{L}_\lambda$ such that $\mathcal{L}_\lambda(\V)=\lambda\V$. Moreover, any two nonzero null vectors are related
by a Lorentz isometry.  Hence, when $\V$ is null, no homothety is needed and the choice of the null vector $\V$ is irrelevant.
\end{remark}

Our first result gives an explicit formula for the sectional curvature of nondegenerate planes in $\R^3_1$.

\begin{proposition}\label{pr23}
Let $\Pi$ be a nondegenerate plane of $\R_1^3$. Then its ambient sectional curvature $K(\Pi)$ with respect to the css-connection determined by $\V$ is
\begin{equation}\label{esp}
K(\Pi)=\frac{1+\alpha^2}{2}(\delta-\varepsilon\langle \V,N\rangle^2)-\alpha\delta,
\end{equation}
where $N$ is the unit normal to $\Pi$ with $\varepsilon=\langle N,N\rangle$ and $\delta=\langle \V,\V\rangle$.
\end{proposition}
\begin{proof}
Let $\{e_1,e_2\}$ be a pseudo-orthonormal basis of $\Pi$, with 
$$\langle e_i,e_j\rangle=\varepsilon_i\delta_{ij}, \qquad \varepsilon_i\in\{-1,1\}.$$ 
Setting $X=e_1$, $Y=Z=e_2$ in \eqref{ss}, and using $\langle e_1,e_2\rangle=0$,
$$
R(e_1,e_2)e_2=\bigl(\eta(e_2)^2-\alpha\delta\,\varepsilon_2\bigr)e_1
-\eta(e_1)\eta(e_2)\,e_2+\alpha^2\eta(e_1)\varepsilon_2\,\V.
$$
Taking the inner product with $e_1$, and using $\langle e_1,e_1\rangle=\varepsilon_1$, $\langle e_2,e_1\rangle=0$, $\langle \V,e_1\rangle=\eta(e_1)$, we obtain
$$
R(e_1,e_2;e_2,e_1)=\varepsilon_1\,\eta(e_2)^2+\alpha^2\varepsilon_2\,\eta(e_1)^2
-\alpha\delta\,\varepsilon_1\varepsilon_2.
$$
Interchanging $e_1$ and $e_2$, that is, taking $X=e_2$, $Y=Z=e_1$, gives
$$
R(e_2,e_1)e_1=\bigl(\eta(e_1)^2-\alpha\delta\,\varepsilon_1\bigr)e_2
-\eta(e_1)\eta(e_2)\,e_1+\alpha^2\eta(e_2)\varepsilon_1\,\V,
$$
and taking the inner product with $e_2$,
\[
R(e_2,e_1;e_1,e_2)=\varepsilon_2\,\eta(e_1)^2+\alpha^2\varepsilon_1\,\eta(e_2)^2
-\alpha\delta\,\varepsilon_1\varepsilon_2.
\]
Dividing each of the two previous identities by $\varepsilon_1\varepsilon_2$, and using $\varepsilon_i^2=1$, we get
\[
\frac{R(e_1,e_2;e_2,e_1)}{\varepsilon_1\varepsilon_2}=\varepsilon_2\eta(e_2)^2+\alpha^2\varepsilon_1\eta(e_1)^2-\alpha\delta
\]
and 
\[
\frac{R(e_2,e_1;e_1,e_2)}{\varepsilon_1\varepsilon_2}=\varepsilon_1\eta(e_1)^2+\alpha^2\varepsilon_2\eta(e_2)^2-\alpha\delta,
\]
respectively. By the definition of $K$ given in \eqref{def1}, we find
$$
K(\Pi)=\frac{1+\alpha^2}{2}(\varepsilon_1\eta(e_1)^2+\varepsilon_2\eta(e_2)^2)-\alpha\delta.
$$
Now, let $\V^\top$ denote the orthogonal projection of $\V$ onto $\Pi$. Then,
$$
\langle \V^\top,\V^\top\rangle=\varepsilon_1\eta(e_1)^2+\varepsilon_2\eta(e_2)^2.
$$
Moreover, since $\V^\top=\V-\varepsilon\eta(N)N$, we have
$$
\langle \V^\top,\V^\top\rangle=\delta-\varepsilon \langle \V,N\rangle^2.
$$
Substituting this identity into the previous expression for $K(\Pi)$ gives the result. 
\end{proof}

It follows from Proposition~\ref{pr23} that the sectional curvature $K(\Pi)$ depends only on $\langle\V^\top,\V^\top\rangle$ and  is therefore independent of the choice of orthonormal basis of $\Pi$.

\begin{remark}\label{rem-null}
Assume that $\V$ is null, so that $\delta=0$. Then \eqref{esp} reduces to $K(\Pi)=-\varepsilon\frac{1+\alpha^2}{2}\langle\V,N\rangle^2$. If $\Pi$ is spacelike, then $\varepsilon=-1$. Since a nonzero null vector cannot be tangent
to a spacelike plane, we have $K(\Pi)>0$. If $\Pi$ is timelike, then $\varepsilon=1$,
so $K(\Pi)\leq0$, with equality if and only if
$\V$ is tangent to $\Pi$.
\end{remark}

We obtain the following consequences of Proposition~\ref{pr23}.

\begin{corollary}\label{cpalpha0}
Let $\nabla$ be the non-metric css-connection determined by $\V$. Then the ambient sectional curvature $K(\Pi)$ of a nondegenerate plane $\Pi\subset \R_1^3$ with respect to $\nabla$ satisfies:
\begin{enumerate}
\item[(i)] If $\Pi$ is spacelike, then
$$
K(\Pi)\geq 0, \quad K(\Pi)> 0, \quad K(\Pi)\geq \frac12,
$$
according as $\V$ is timelike, null, or spacelike, respectively.
\item[(ii)] If $\Pi$ is timelike, then
$$
K(\Pi)\leq -\frac12, \quad K(\Pi)\leq 0, \quad K(\Pi)\leq \frac12,
$$
according as $\V$ is timelike, null, or spacelike, respectively.
\end{enumerate}
\end{corollary}
\begin{proof}
By  Proposition \ref{pr23}, setting $\alpha=0$ yields
$$
K(\Pi)=\frac12(\delta-\varepsilon\langle\V,N\rangle^2).
$$ 
\begin{enumerate}
\item[(i)] Assume that $\Pi$ is spacelike, i.e.,  $\varepsilon=-1$. Then, $K(\Pi)=\frac12(\delta+\langle\V,N\rangle^2)$. If $\V$ is timelike, then $\delta=-1$. Since both $\V$ and $N$ are timelike, it follows that $\langle\V,N\rangle^2\geq 1$. Hence, $K(\Pi)=\frac12(-1+\langle\V,N\rangle^2)\geq 0$. If $\V$ is null ($\delta=0$), then $\V$ cannot lie in $\Pi$ and hence $\langle\V,N\rangle^2>0$, implying $K(\Pi)>0$. If $\V$ is spacelike ($\delta=1$), then $K(\Pi)=\frac12(1+\langle\V,N\rangle^2)\geq \frac12$.
\item[(ii)] Assume that $\Pi$ is timelike, i.e.,  $\varepsilon=1$. Then, $K(\Pi)=\frac12(\delta-\langle\V,N\rangle^2)$. If $\V$ is timelike ($\delta=-1$), then $K(\Pi)=-\frac12(1+\langle\V,N\rangle^2)\leq -\frac12$.  If $\V$ is null ($\delta=0$), then $K(\Pi)=-\frac12\langle\V,N\rangle^2\leq 0$. Finally, if $\V$ is spacelike ($\delta=1$), then $K(\Pi)=\frac12(1-\langle\V,N\rangle^2)\leq \frac12$.
\end{enumerate}
\end{proof}

\begin{corollary}\label{cpalpha1}
Let $\nabla$ be the metric css-connection determined by $\V$. Then the ambient sectional curvature $K(\Pi)$ of a nondegenerate plane $\Pi\subset \R_1^3$ with respect to $\nabla$ satisfies:
\begin{enumerate}
\item[(i)] If $\Pi$ is spacelike, then
$$
K(\Pi)\geq 1, \quad K(\Pi)> 0, \quad K(\Pi)\geq 0,
$$
according as $\V$ is timelike, null, or spacelike, respectively.
\item[(ii)] If $\Pi$ is timelike, then $K(\Pi)\leq 0$.
\end{enumerate}
\end{corollary}
\begin{proof}
Writing $\alpha=1$ in \eqref{esp} yields $K(\Pi)=-\varepsilon\langle \V,N\rangle^2$.
\begin{enumerate}
\item[(i)] Assume that $\Pi$ is spacelike, i.e.,  $\varepsilon=-1$. Then, $K(\Pi)=\langle\V,N\rangle^2$. If $\V$ is timelike, then $K(\Pi)=\langle\V,N\rangle^2\geq 1$. If $\V$ is null, then $K(\Pi)>0$, and if $\V$ is spacelike, then $K(\Pi)\geq0$. 
\item[(ii)] Assume that $\Pi$ is timelike, i.e.,  $\varepsilon=1$. Then, $K(\Pi)=-\langle \V,N\rangle^2$. Since $\langle \V,N\rangle^2\geq0$, it follows that  $K(\Pi)\leq 0$.
\end{enumerate}
\end{proof}

We finish this section by determining the scalar curvature of $\R_1^3$.

\begin{corollary}\label{cscalar}
Let $\nabla$ be the css-connection on $\R_1^3$ determined by $\V$ and $\delta=\langle \V,\V\rangle.$ Then, the scalar curvature $S$ of $\R_1^3$ with respect to $\nabla$ is
$$
S=
\begin{cases}
\delta, & \alpha=0,\\ 
-\delta, & \alpha=1.
\end{cases}
$$
In particular, $S=0$ when $\V$ is null.
\end{corollary}

\begin{proof}
Let $\{e_1,e_2,e_3\}$ be a pseudo-orthonormal basis of $\R_1^3$, where 
$$
\langle e_i,e_j\rangle=\varepsilon_i\delta_{ij},
\qquad
\varepsilon_i\in\{-1,1\}.
$$
For $1\leq i<j\leq3$, set
$\Pi_{ij}=\operatorname{span}\{e_i,e_j\}$.
We adopt the convention $S=\sum_{i<j}K(\Pi_{ij})$, which is half the scalar curvature defined
by the metric trace of the Ricci tensor. For each pair $i<j$, let $e_k$ be the unit normal to $\Pi_{ij}$. By Proposition~\ref{pr23},
$$
K(\Pi_{ij})
=
\frac{1+\alpha^2}{2}
\left(
\delta-\varepsilon_k\langle\V,e_k\rangle^2
\right)
-\alpha\delta.
$$
Summing over all pairs $i<j$, we obtain
\begin{align*}
S
&=
\frac{1+\alpha^2}{2}
\left(
3\delta-
\sum_{k=1}^{3}
\varepsilon_k\langle\V,e_k\rangle^2
\right)
-3\alpha\delta,
\end{align*}
or equivalently,
\begin{align*}
S
&=
\frac{1+\alpha^2}{2}(3\delta-\delta)
-3\alpha\delta\\
&=
(1+\alpha^2-3\alpha)\delta.
\end{align*}
For $\alpha=0$ and $\alpha=1$, respectively, this yields $S=\delta$ and $S=-\delta.$ This completes the proof.
\end{proof}

\section{Curvature of nondegenerate surfaces with respect to css-connections} \label{sec3}

Let $\Sigma\subset\R^3_1$ be a nondegenerate immersed surface. Let $\nabla^0$ denote the Levi-Civita connection of $\R^3_1$ and
$\nabla^{0_\Sigma}$ the induced Levi-Civita connection on $\Sigma$. The Gauss formula is
given by
\begin{equation}\label{gauss0}
\nabla^0_XY=\nabla^{0_\Sigma}_XY+h^0(X,Y), \qquad X,Y\in \X(\Sigma),
\end{equation}
where $h^0$ is the second fundamental form of the immersion with respect to $\nabla^0$. Since $\Sigma$ is nondegenerate, we can define a unit normal vector field $N$ on $\Sigma$; it is spacelike if $\Sigma$ is timelike, and timelike if $\Sigma$ is spacelike. 

Let $\varepsilon=\langle N,N\rangle$. The mean curvature $H^0$ of $\Sigma$ is defined by 
$$
H^0=\frac{\varepsilon}{2}\langle\operatorname{trace}(h^0),N\rangle.
$$ 
The shape operator $A_p^0$, $p\in\Sigma$, is the self-adjoint endomorphism of $T_p\Sigma$ defined by 
$$\langle A_p^0(X_p),Y_p\rangle
=\langle h_p^0(X_p,Y_p),N(p)\rangle,
\qquad X_p,Y_p\in T_p\Sigma.
$$
Consequently, 
$$H^0=\frac{\varepsilon}{2}\operatorname{trace}(A^0),
\qquad
G=\varepsilon\det(A^0),
$$
where $G$ denotes the Gauss curvature of $\Sigma$. If $\varepsilon=-1$, then the shape operator $A^0$ is diagonalizable and its eigenvalues are called the principal curvatures of $\Sigma$.

Let $\Psi(u,v)$ be a local parametrization of $\Sigma$. The coefficients of the first fundamental form are
\[
e=\langle \Psi_u,\Psi_u\rangle,\qquad
f=\langle \Psi_u,\Psi_v\rangle,\qquad
g=\langle \Psi_v,\Psi_v\rangle.
\]
The coefficients of the second fundamental form $h^0$ are
\[
l^0=\langle N,\Psi_{uu}\rangle,\qquad
m^0=\langle N,\Psi_{uv}\rangle,\qquad
n^0=\langle N,\Psi_{vv}\rangle.
\]
The mean curvature $H^0$ and the Gauss curvature $G$ are given by
\[
H^0=\frac{\varepsilon}{2}\,
\frac{l^0g-2m^0f+n^0e}{eg-f^2},
\qquad
G=\varepsilon\,
\frac{l^0n^0-(m^0)^2}{eg-f^2}.
\]

Let $\nabla$ be the css-connection on $\R^3_1$ defined by \eqref{con0}. The Gauss formula for the surface $\Sigma\subset\R_1^3$ with respect to $\nabla$ is
\begin{equation}\label{gauss1}
\nabla_XY=\nabla_X^{\Sigma}Y+h(X,Y),
\qquad X,Y\in\X(\Sigma),
\end{equation}
where $\nabla^{\Sigma}$ denotes the connection on $\Sigma$ induced by $\nabla$ and $h$ is the second fundamental form of $\Sigma$ with respect to $\nabla$.

Set $\V=\V^\top+\V^\perp$, where $\V^\top$ and $\V^\perp$ denote the tangential and normal components of $\V$ along $\Sigma$, respectively. Since
$\V^\perp=\varepsilon\eta(N)N$,
comparing the Gauss formulas \eqref{gauss0}
and \eqref{gauss1}, together with \eqref{con0}, yields
$$
 h(X,Y)= h^0(X,Y)-\alpha\varepsilon\eta(N)\langle X,Y\rangle N, \qquad X,Y\in \X(\Sigma),
$$
or equivalently,
$$
l=l^0-\alpha\eta(N)e, \quad m=m^0-\alpha\eta(N)f, \quad n=n^0-\alpha\eta(N)g, 
$$
where $l,m,n$ and $l^0,m^0,n^0$ denote the coefficients of the second fundamental forms $h$ and $h^0$, respectively.

The next result expresses the intrinsic sectional curvature $K^\Sigma$ of $\Sigma\subset\R_1^3$ with respect to the connection \eqref{con0} in terms of its Gauss curvature $G$ and mean curvature $H^0$.

\begin{proposition}\label{pr31}
Let $\Sigma$ be a nondegenerate surface in
$\R_1^3$ with unit normal $N$ and let
$\varepsilon=\langle N,N\rangle$. Then the intrinsic sectional curvature $K^\Sigma$ of $\Sigma $ with respect to the  css-connection $\nabla$ determined by $\V$ is 
\begin{equation} \label{sectcurvat}
\begin{aligned}
K^\Sigma&=K+G-(1+\alpha)\eta(N)H^0 +\varepsilon\alpha\eta(N)^2\\
&=\frac{(1-\alpha)^2}{2}(\delta-\varepsilon\eta(N)^2)+G-(1+\alpha)\eta(N)H^0,
\end{aligned}
\end{equation}
where $K$ is the ambient sectional curvature of $\Sigma$ with respect to $\nabla$,  $G$ and $H^0$ are, respectively,  the Gauss and mean curvatures of $\Sigma$ with respect to the Levi-Civita connection, $\eta(N)=\langle \V,N\rangle$ and $\delta=\langle \V,\V\rangle$.
\end{proposition}

\begin{proof}
Let $X,Y,Z\in \X(\Sigma)$. Denote by $R$ and
$R^\Sigma$ the curvature tensors of $\R_1^3$ and $\Sigma$ associated
with the css-connection $\nabla$ and the induced connection
$\nabla^\Sigma$, respectively. Since $\langle X,N\rangle=0$, the Weingarten
formula $\nabla_X^0N=-A^0X$, together with
\eqref{con0}, yields $\nabla_XN=-\bigl(A^0X-\eta(N)X\bigr)$. Taking tangential components in the
curvature identity and using the Gauss
formula \eqref{gauss1}, we obtain
$$
\begin{aligned}
R^\Sigma(X,Y)Z
={}&\bigl(R(X,Y)Z\bigr)^\top+\varepsilon\langle h(Y,Z),N\rangle
 \bigl(A^0X-\eta(N)X\bigr)\\
&-\varepsilon\langle h(X,Z),N\rangle
 \bigl(A^0Y-\eta(N)Y\bigr).
\end{aligned}
$$
Let $\{e_1,e_2\}$ be a local pseudo-orthonormal tangent frame on $\Sigma$ satisfying
$$
\langle e_i,e_j\rangle=\varepsilon_i\delta_{ij},
\qquad
\varepsilon_i\in\{-1,1\},
\qquad
\varepsilon_1\varepsilon_2=-\varepsilon.
$$
Setting $X=e_1$, $Y=e_2$, and $Z=e_2$, and then taking the inner
product with $e_1$, we obtain
$$
\begin{aligned}
R^\Sigma(e_1,e_2;e_2,e_1)
=R(e_1,e_2;e_2,e_1)+\varepsilon
\bigl(n^0-\alpha\varepsilon_2\eta(N)\bigr)
\bigl(l^0-\varepsilon_1\eta(N)\bigr)
-\varepsilon(m^0)^2.
\end{aligned}
$$
Similarly, 
$$
\begin{aligned}
R^\Sigma(e_2,e_1;e_1,e_2)
=R(e_2,e_1;e_1,e_2)+\varepsilon
\bigl(l^0-\alpha\varepsilon_1\eta(N)\bigr)
\bigl(n^0-\varepsilon_2\eta(N)\bigr)
-\varepsilon(m^0)^2.
\end{aligned}
$$
Here $l^0=\langle A^0e_1,e_1\rangle$, $m^0=\langle A^0e_1,e_2\rangle$ and $n^0=\langle A^0e_2,e_2\rangle$. Adding these two identities, dividing by
$2\varepsilon_1\varepsilon_2$ and using
$G=\varepsilon\varepsilon_1\varepsilon_2\bigl(l^0n^0-(m^0)^2\bigr)$ and $\operatorname{trace}(A^0)=\varepsilon_1l^0+\varepsilon_2n^0=2\varepsilon H^0$, we obtain
$$
K^\Sigma=K+G-(1+\alpha)\eta(N)H^0+\varepsilon\alpha\eta(N)^2.
$$
 Here we have also used   the fact that the terms linear in $\eta(N)$ involve $\varepsilon_1n^0+\varepsilon_2l^0$, which is equal to $-\varepsilon\operatorname{trace}(A^0)=-2H^0$ because $\varepsilon_1\varepsilon_2=-\varepsilon$. The second expression in \eqref{sectcurvat} is a direct consequence of \eqref{esp}.
\end{proof}

The formula \eqref{sectcurvat} can also be written in intrinsic terms. Let $\Psi$ denote the position vector of $\Sigma$ and let $\xi:\Sigma\to\R$ be the height function defined by $\xi=\langle\V,\Psi\rangle$. Since $\V$ is constant, the gradient of $\xi$ with respect to the induced metric is $\nabla^{0_\Sigma}\xi=\V^\top$.

\begin{corollary}\label{intrinsic}
With the notation of Proposition \ref{pr31}, the intrinsic sectional curvature of $\Sigma$ with respect to the css-connection is
\begin{equation}\label{intr}
K^\Sigma=G-\frac{1+\alpha}{2}\,\Delta^\Sigma \xi+\frac{(1-\alpha)^2}{2}\,\langle\nabla^{0_\Sigma}\xi,\nabla^{0_\Sigma}\xi\rangle,
\end{equation}
where $\Delta^\Sigma$ is the Laplacian of $\Sigma$. Hence, $K^\Sigma$ depends only on the induced metric of $\Sigma$ and on the function $\xi$. 
\end{corollary}

\begin{proof}
First, since $\nabla^{0_\Sigma}\xi=\V^\top$, we have $\langle\nabla^{0_\Sigma}\xi,\nabla^{0_\Sigma}\xi\rangle=\delta-\varepsilon\eta(N)^2. $ Let $\{e_1,e_2\}$ be a pseudo-orthonormal tangent frame with $\varepsilon_i=\langle e_i,e_i\rangle$. Since $\V$ is constant, using $\V^\top=\V-\varepsilon\eta(N)N$ and the Weingarten formula $\nabla^0_XN=-A^0X$, we obtain
$$
\Delta^\Sigma \xi=\operatorname{div_\Sigma}(\V^\top)=-\varepsilon\eta(N)\sum_{i=1}^2\varepsilon_i\langle\nabla^0_{e_i}N,e_i\rangle=\varepsilon\eta(N)\operatorname{trace}(A^0)=2\eta(N)H^0.
$$
Substituting these identities into the second expression of \eqref{sectcurvat} gives \eqref{intr}.
\end{proof}

In the metric case ($\alpha=1$), \eqref{intr} reduces to $K^\Sigma=G-\Delta^\Sigma \xi$, while in the non-metric case ($\alpha=0$) it can be written as $K^\Sigma=G+\frac12\,e^{\xi}\Delta^\Sigma(e^{-\xi})$.

As immediate consequences of Proposition~\ref{pr31} and formula \eqref{sectcurvat}, we obtain the following particular cases.
\begin{corollary}\label{cor32}
Let $\Sigma$ be a nondegenerate surface of $\R_1^3$ and let $K^\Sigma$ and $K$ denote the intrinsic and the ambient sectional curvatures of $\Sigma$ with respect to the css-connection determined by $\V$, respectively. Then:
\begin{enumerate}
\item[(i)] If $\Sigma$ is a plane, then $K^\Sigma=K+\varepsilon\alpha\eta(N)^2=\frac{(1-\alpha)^2}{2}\big(\delta-\varepsilon\eta(N)^2\big)$, which is constant.  In particular, $K^\Sigma=K$ in the non-metric case, while $K^\Sigma=0$ in the metric case.
\item[(ii)] If $\V$ is tangent to $\Sigma$, then $\Sigma$ is an open piece of a cylindrical surface whose rulings are parallel to $\V$. Moreover, $G=0$ and
$$
K^\Sigma=K=\delta\,\frac{(1-\alpha)^2}{2}.
$$
This holds for every causal character of $\V$, and $K^\Sigma=0$ when $\V$ is null.
\item[(iii)] If $\Sigma$ is minimal, i.e.,  $H^0=0$, then $K^\Sigma=K+G+\varepsilon\alpha\eta(N)^2$.
\end{enumerate}
\end{corollary}

\begin{proof}
\begin{enumerate}
\item[(i)] If $\Sigma$ is a plane, it is totally geodesic with respect to the ambient
Levi-Civita connection. Thus, $h^0\equiv0$ and $A^0=0$, implying $G=\det(A^0)=0$ and $H^0=\frac{\varepsilon}{2}\operatorname{trace}(A^0)=0$. Substituting these identities into
\eqref{sectcurvat}, we obtain $K^\Sigma=K+\varepsilon\alpha\eta(N)^2$. It follows from \eqref{esp} that $K^\Sigma=\frac{(1-\alpha)^2}{2}\big(\delta-\varepsilon\eta(N)^2\big)$.  Since $N$ is constant, $K^\Sigma$ is
constant. In particular, $K^\Sigma=K$
when $\alpha=0$, while $K^\Sigma=0$
when $\alpha=1$.

\item[(ii)] Assume $\eta(N)=\langle\V,N\rangle=0$ on $\Sigma$. Since $\V$ is a constant vector field tangent to $\Sigma$, its integral curves in $\Sigma$ are segments of the straight lines $q+t\V$, $q\in\Sigma$. Hence, $\Sigma$ is locally parametrized by $\psi(s,t)=\gamma(s)+t\V$, where $\gamma$ is a curve in $\Sigma$ transversal to $\V$. Since $\psi_t=\V$ and $\psi_{st}=\psi_{tt}=0$, it follows that $m^0=n^0=0$ and $G=\varepsilon\,(l^0n^0-(m^0)^2)/(eg-f^2)=0$. As $\eta(N)=0$, formula \eqref{sectcurvat} yields $K^\Sigma=K$. Furthermore,
\eqref{esp} gives $K=\frac{1+\alpha^2}{2}\delta-\alpha\delta=\delta\frac{(1-\alpha)^2}{2}$.

\item[(iii)] This follows immediately from
\eqref{sectcurvat} by setting
$H^0=0$.
\end{enumerate}
\end{proof}

We next give an explicit family of translation surfaces satisfying $K^\Sigma=K$, providing a Lorentzian counterpart of Example 2.6 in \cite{alm1}.

\begin{example}\label{translation}
Let $\V=\partial_z$, so that $\delta=\langle\V,\V\rangle=-1$, and let $\alpha=0$. Consider the nondegenerate translation surfaces $\Sigma$ given by $z=u(x,y)=f(x)+g(y)$, where $f=f(x)$, $g=g(y)$ are smooth functions. We look for such surfaces satisfying $K^\Sigma=K$. By \eqref{sectcurvat}, the condition $K^\Sigma=K$ is equivalent to
\begin{equation}\label{gnh}
G-\eta(N)H^0=0.
\end{equation}
\begin{enumerate}
\item[(i)] The case where $\Sigma$ is spacelike. Let $W:=1-u_x^2-u_y^2>0$ and let $N$ be the unit normal
$$
N=\frac{1}{\sqrt W}(u_x,u_y,1),\qquad \varepsilon=\langle N,N\rangle=-1.
$$
Using $u_{xy}=0$, one obtains
\begin{equation*}
\begin{split}
\eta(N)&=-\frac{1}{\sqrt W},\\
G&=-\frac{u_{xx}u_{yy}}{W^2},\\
H^0&=\frac{u_{xx}(1-u_y^2)+u_{yy}(1-u_x^2)}{2W^{3/2}}.
\end{split}
\end{equation*}
Substituting these into \eqref{gnh} and multiplying by $W^2$, we obtain
\[
u_{xx}u_{yy}=\frac12\left(u_{xx}(1-u_y^2)+u_{yy}(1-u_x^2)\right),
\]
that is,
 \begin{equation}\label{tra}
2f''g''=f''(1-g'^2)+g''(1-f'^2).
\end{equation}
One solution of \eqref{tra} is obtained  when $f''=g''\equiv0$. Then $\Sigma$ is a plane, in agreement with  item (i) of Corollary~\ref{cor32}. If $g''=0$ and $f''\neq0$ on some interval, then \eqref{tra} gives $g'^2=1$ and so $W=-f'^2\leq0$, contradicting $W>0$. The case $f''=0$ and $g''\neq0$ is analogous. We discard these cases and assume $f''g''\neq0$ on $I\times J$. Dividing \eqref{tra} by $f''g''$, we obtain
$$
\frac{1-g'^2}{g''}+\frac{1-f'^2}{f''}=2. 
$$
Separating variables, there is a nonzero constant $c$ such that
$$
\frac{1-f'^2}{f''}=\frac1c,\qquad \frac{1-g'^2}{g''}=2-\frac1c=\frac{2c-1}{c}.
$$
Here $c\neq\frac12$, since otherwise $g'^2=1$ and thus $g''=0$, contradicting $f''g''\neq0$. Hence
$$
f''=c(1-f'^2),\qquad g''=k(1-g'^2),\qquad k=\frac{c}{2c-1}.
$$
Since $W>0$ implies $|f'|<1$ and $|g'|<1$, the solutions are $f'=\tanh(c(x-x_0))$ and $g'=\tanh(k(y-y_0))$. Integrating the two separated ODEs, we obtain, up to translations of $x,y$ and additive constants,
\begin{equation}\label{solex}
\begin{split}
f(x)&=\frac1c\log\cosh(cx),\\
g(y)&=\frac{1}{k}\log\cosh(ky),\qquad k=\frac{c}{2c-1},\quad c\neq0,\tfrac12.
\end{split}
\end{equation}
Since $f'(0)=g'(0)=0$, we have $W>0$ near the origin. In contrast with the Euclidean solutions in \cite{alm1}, the functions $f,g$ in \eqref{solex} have no singularities.
\item[(ii)] The case where $\Sigma$ is timelike. Set $W:=u_x^2+u_y^2-1>0$. Choosing $N=W^{-1/2}(u_x,u_y,1)$, we have $\varepsilon=1$ and the determinant of the first fundamental form is $-W$. A direct computation shows that $\eta(N)$, $G$ and $H^0$ are given by the same expressions as in the spacelike case. Hence we arrive again at \eqref{tra}.

First, assume that $f''g''\neq0$. Separating variables as above, we obtain $f''=c(1-f'^2)$ and $g''=k(1-g'^2)$ with $k=\frac{c}{2c-1}$, where $c\neq\frac12$ as in the spacelike case. Unlike the spacelike case, the timelike condition allows $|f'|>1$ or $|g'|>1$. Hence, besides $f'=\tanh(cx)$, we have the branch $f'=\coth(cx)$, and similarly for $g$. Up to translations of $x,y$ and additive constants, $f$ is one of
$$
\frac1c\log\cosh(cx),\qquad \frac1c\log|\sinh(cx)|,
$$
and $g$ is one of
$$
\frac{1}{k}\log\cosh(ky),\qquad \frac{1}{k}\log|\sinh(ky)|.
$$
Each combination defines a timelike translation surface on any nonempty open region where $f'(x)^2+g'(y)^2>1$.

Assume now that $f''g''=0$ on some open set, say $g''=0$ and $f''\neq0$. Then \eqref{tra} gives $g'^2=1$, that is, $g(y)=\pm y$ up to an additive constant, and  \eqref{tra} then holds for any function $f$. For $u(x,y)=f(x)\pm y$ we have $u_{yy}=0$ and $1-u_y^2=0$, so that $G=H^0=0$. These surfaces are defined where $f'\neq0$, so that $W=f'^2>0$. They are cylindrical surfaces whose rulings are parallel to the null vector $(0,1,\pm1)$ (see Proposition~\ref{nullrulings}), and they have no spacelike counterpart. Finally, the case $c=\frac12$ does not yield an isolated surface in this family, with $f$ arbitrary.
\end{enumerate}
\end{example}

\section{Cylindrical surfaces of constant sectional curvature with respect to css-connections} \label{sec4}

We consider two classes of surfaces in $\R_1^3$ invariant under a one-parameter group of isometries: cylindrical surfaces, which are invariant under translations, and rotational surfaces, which are invariant under rotations.
We study the former in this section and the latter in Section~\ref{sec5}.

Let $w\in\R_1^3$ be a fixed unit vector, with $\varepsilon_2=\langle w,w\rangle\in\{-1,1\}$, and let $\gamma:I\to\R_1^3$ be a curve parametrized by arc length contained in a plane orthogonal to $w$, that is, $\langle\gamma',w\rangle=0$ and $\langle\gamma',\gamma'\rangle=\varepsilon_1\in\{-1,1\}$ identically. Let $\mathbf n$ be the unit vector field orthogonal to both $\gamma'$ and $w$, so that $\gamma''=\kappa\,\mathbf n$ for some function $\kappa$, and set $\varepsilon_n=\langle\mathbf n,\mathbf n\rangle$. Then, $\varepsilon_1\varepsilon_2\varepsilon_n=-1$. Consequently, a cylindrical surface with rulings parallel to $w$ and generating curve $\gamma$ is parametrized by
\begin{equation}\label{eqcyl}
\begin{cases} 
      \psi(s,t)=\gamma(s)+t\,w,\quad s\in I,\ t\in\R. \\ 
      \langle\gamma',\gamma'\rangle=\varepsilon_1, \quad \langle w,w\rangle=\varepsilon_2, \quad \langle\mathbf n,\mathbf n\rangle=\varepsilon_n, \quad \varepsilon_1,\varepsilon_2,\varepsilon_n\in\{-1,1\}.
\end{cases}
\end{equation}
The unit normal of $\Sigma$ is $N=\mathbf n$.  

Cylindrical surfaces with null rulings are considered separately at the end of this section, in Proposition~\ref{nullrulings}.

\begin{proposition}\label{thcyl}
Let $\Sigma$ be a cylindrical surface of $\R^3_1$ parametrized by \eqref{eqcyl}. Then its intrinsic sectional curvature $K^\Sigma$ with respect to the css-connection determined by $\V$ is
\begin{equation}\label{15}
K^\Sigma=\frac{1+\alpha^2}{2}\Big(\varepsilon_1\,\eta(\gamma')^2+\varepsilon_2\,\eta(w)^2\Big)-\alpha\delta
-(1+\alpha)\,\eta(\mathbf n)\,\frac{\kappa}{2\varepsilon_1}
+\varepsilon_n\alpha\,\eta(\mathbf n)^2,
\end{equation}
where $\eta(X)=\langle\V,X\rangle$.
\end{proposition}

\begin{proof}
Since $\psi_t=w$ is constant, $\psi_{tt}=0$ and $\psi_{st}=0$. Also, $\psi_{ss}=\gamma''=\kappa\,\mathbf n$. Hence the coefficients of the first and second fundamental forms of $\Sigma$ are
\[
e=\varepsilon_1,\quad f=0,\quad g=\varepsilon_2,\]
\[l^0=\kappa\varepsilon_n,\quad m^0=0,\quad n^0=0.
\]
Consequently, the Gauss and mean curvatures
with respect to the Levi-Civita connection are
\[
G=0,\qquad
H^0=\frac{\varepsilon_n}{2}
\frac{l^0g}{eg-f^2}
=\frac{\kappa}{2\varepsilon_1}.
\]
On the other hand, since $\{\gamma',w\}$ is a pseudo-orthonormal basis of the tangent plane of $\Sigma$, Proposition  \ref{pr23} gives
$$
K=\frac{1+\alpha^2}{2}\big(\varepsilon_1\eta(\gamma')^2+\varepsilon_2\eta(w)^2\big)-\alpha\delta.
$$
Substituting the expressions for $G$, $H^0$ and $K$ into \eqref{sectcurvat}, with $\eta(N)=\eta(\mathbf n)$ and $\varepsilon=\varepsilon_n$, yields \eqref{15}.
\end{proof}

The expression \eqref{15} can be simplified by introducing the functions
\begin{equation}\label{int-func-pc}
p:=\eta(\gamma')=\langle\V,\gamma'\rangle,\qquad c:=\eta(w)=\langle\V,w\rangle.
\end{equation}

\begin{corollary}\label{thp}
Let $\Sigma$ be a cylindrical surface of $\R_1^3$ parametrized by \eqref{eqcyl}. Then its intrinsic sectional curvature with respect to the css-connection is given by
\begin{equation}\label{Kp}
K^\Sigma=\frac{(1-\alpha)^2}{2}\big(\varepsilon_1p^2+\varepsilon_2c^2\big)-\frac{1+\alpha}{2}\,\varepsilon_1\,p'.
\end{equation}
\end{corollary}

\begin{proof}
Since $\{\gamma',w,\mathbf n\}$ is a pseudo-orthonormal frame, we have $\V=\varepsilon_1p\,\gamma'+\varepsilon_2c\,w+\varepsilon_n\eta(\mathbf n)\mathbf n$. Hence, $\delta=\varepsilon_1p^2+\varepsilon_2c^2+\varepsilon_n\eta(\mathbf n)^2$. Moreover, $p'=\langle\V,\gamma''\rangle=\kappa\,\eta(\mathbf n)$. Substituting both identities into \eqref{15},
$$
K^\Sigma=\frac{1+\alpha^2}{2}\big(\varepsilon_1p^2+\varepsilon_2c^2\big)-\alpha\big(\varepsilon_1p^2+\varepsilon_2c^2\big)-\frac{1+\alpha}{2\varepsilon_1}\,p',
$$
which reduces to \eqref{Kp} because $1/\varepsilon_1=\varepsilon_1$.
\end{proof}

Corollary~\ref{thp} allows us to classify all cylindrical surfaces with constant intrinsic sectional curvature, for any position of the rulings with respect to $\V$. By Corollary~\ref{cor32}(i), planes have constant intrinsic sectional curvature.

\begin{theorem}\label{clascyl}
Let $\Sigma$ be a non-planar cylindrical surface of $\R^3_1$ parametrized by \eqref{eqcyl} and let $K^\Sigma$ denote the intrinsic sectional curvature with respect to the css-connection determined by $\V$. Then $K^\Sigma$ is constant, $K^\Sigma=K_0$, if and only if one of the following holds:
\begin{enumerate}
\item[(i)] $\V$ is parallel to $w$. In this case, $K_0=\delta\frac{(1-\alpha)^2}{2}$ and $\gamma$ is arbitrary.
\item[(ii)] The connection is non-metric ($\alpha=0$)  and, up to a translation of $s$, the function $p$ defined by \eqref{int-func-pc} is one of the following, where $\mu:=\varepsilon_1\big(2K_0-\varepsilon_2c^2\big)$:
\begin{enumerate}
\item[(ii.1)]  if $\mu>0$, with $a=\sqrt{\mu}$: $p(s)=-a\tanh(as)$ or $p(s)=-a\coth(as)$.
\item[(ii.2)]  if $\mu=0$: $p(s)=-1/s$.
\item[(ii.3)]  if $\mu<0$, with $b=\sqrt{-\mu}$: $p(s)=b\tan(bs)$.
\end{enumerate}
\item[(iii)]  The connection is metric ($\alpha=1$), $K_0\neq0$ and, up to a translation of $s$, $p(s)=-\varepsilon_1K_0\,s$.
\end{enumerate}
Moreover, in cases (ii) and (iii) the curve $\gamma$ is determined by $p$ up to a translation, as follows. Let $P=w^\perp$, which contains $\gamma'$, let $\V_P=\V-\varepsilon_2c\,w$ be the orthogonal projection of $\V$ onto $P$, and $\delta_P=\langle\V_P,\V_P\rangle$.
\begin{enumerate}
\item[(I)] If $\delta_P\neq0$, let $u\in P$ be a unit vector orthogonal to $\V_P$ and $\varepsilon_u=\langle u,u\rangle$. Then
$$
\gamma'=\frac{p}{\delta_P}\,\V_P+\beta\,u,\qquad \beta=\pm\sqrt{\varepsilon_u\Big(\varepsilon_1-\frac{p^2}{\delta_P}\Big)},
$$
and $\gamma$ is defined on the intervals where $\varepsilon_u\big(\varepsilon_1-p^2/\delta_P\big)>0$.
\item[(II)] If $\delta_P=0$ (so $\V_P\neq0$ is null and $P$ is timelike), let $L\in P$ be the null vector with $\langle\V_P,L\rangle=1$. Then
$$
\gamma'=\frac{\varepsilon_1}{2p}\,\V_P+p\,L,
$$
and $\gamma$ is defined on the intervals where $p\neq0$.
\end{enumerate}
\end{theorem}

\begin{proof}
If $\V$ is parallel to $w$, then $p=0$ and \eqref{Kp} gives $K^\Sigma=\frac{(1-\alpha)^2}{2}\varepsilon_2c^2=\delta\frac{(1-\alpha)^2}{2}$, because $\V=\varepsilon_2c\,w$ and $\delta=\varepsilon_2c^2$. This is case (i). Assume now that $\V$ is not parallel to $w$. By \eqref{Kp}, $K^\Sigma=K_0$ if and only if
\begin{equation}\label{ricc}
p'=p^2-\mu\quad (\alpha=0),\qquad p'=-\varepsilon_1K_0\quad (\alpha=1).
\end{equation}
For $\alpha=0$, the non-constant solutions of the Riccati equation $p'=p^2-\mu$ are, up to a translation of $s$, those listed in (ii). The constant solutions are $p=\pm\sqrt\mu$ if $\mu\geq0$. For $\alpha=1$, $p$ is non-constant if and only if $K_0\neq0$, and then $p=-\varepsilon_1K_0(s-s_0)$, which proves case (iii).

Since $\V$ is not parallel to $w$, it follows that $\V_P\neq0$. Moreover, since $\gamma'\in P$ and $\langle\V_P,\gamma'\rangle=\langle\V,\gamma'\rangle=p$, we distinguish two cases. If $\delta_P\neq0$, then $\{\V_P,u\}$ is an orthogonal basis of $P$, so $\gamma'=\frac{p}{\delta_P}\V_P+\beta u$ with $\beta=\varepsilon_u\langle\gamma',u\rangle$. Hence, $\langle\gamma',\gamma'\rangle=\varepsilon_1$ gives $\beta^2=\varepsilon_u(\varepsilon_1-p^2/\delta_P)$. The sign of $\beta$ is constant on each interval where $\beta\neq0$. If $\delta_P=0$, write $\gamma'=a\V_P+bL$. Then $p=\langle\V_P,\gamma'\rangle=b$ and $\varepsilon_1=\langle\gamma',\gamma'\rangle=2ab$, hence $p\neq0$ and $a=\varepsilon_1/(2p)$. This proves  cases (I) and (II). In both cases, if $p$ is constant on an interval, then $\gamma'$ is also constant: in case (I), $\beta^2$ is constant and hence, by continuity, so is $\beta$; in case (II), this is immediate from the expression for $\gamma'$. Hence $\gamma$ is a straight line and $\Sigma$ is a plane. Consequently, $p$ is a non-constant solution of \eqref{ricc}. Conversely, let $p$ be as in (ii) or (iii). The curve $\gamma$ given by (I) or (II) is parametrized by arc length and lies in a plane orthogonal to $w$. It satisfies $\langle\V,\gamma'\rangle=p$, and \eqref{Kp} gives $K^\Sigma=K_0$.
\end{proof}

\begin{remark}
When $\V$ is null, the parametrization \eqref{eqcyl} with $w=\V$ is not applicable, because $\V$ cannot be normalized and the plane $\V^\perp$ is degenerate. However, the conclusion of case (i) of Theorem \ref{clascyl} is still valid for nondegenerate cylindrical surfaces whose rulings are parallel to $\V$. Indeed, by Corollary~\ref{cor32}(ii), such surfaces
satisfy $G=0$ and $K^\Sigma=\delta\,\frac{(1-\alpha)^2}{2}=0.$ See also Proposition~\ref{nullrulings}.
\end{remark}

We next consider cylindrical surfaces whose rulings are orthogonal to $\V$. We distinguish cases according to the causal character of $\V$ and to whether the css-connection is metric or non-metric.

\subsection{The non-metric case ($\alpha=0$)}
We first consider cylindrical surfaces  $\Sigma$ in the non-metric case $\alpha=0$. Since the rulings are assumed to be orthogonal to $\V$, we have $c=\langle\V,w\rangle=0$. When $K^\Sigma$ is a constant $K_0$,  Theorem \ref{clascyl}(ii) gives
\begin{equation}\label{pRicc-solu}
p(s)=
\begin{cases}
-a\tanh(as)\ \text{or}\ -a\coth(as),
& \varepsilon_1K_0>0,\\[1mm]
-1/s,
& K_0=0,\\[1mm]
b\tan(bs),
& \varepsilon_1K_0<0.
\end{cases}
\end{equation}
where $a=\sqrt{2\varepsilon_1K_0}$ and $b=\sqrt{-2\varepsilon_1K_0}$ in the respective cases. Although these solutions are independent of the causal character of $\V$, the resulting generating curves depend,  through the relation $p=\langle\V,\gamma'\rangle$, on whether $\V$ is spacelike, timelike or null.

First, assume that $\V$ is timelike, so $\delta=-1$. Since linear isometries of $\R_1^3$ preserve the metric, the connection \eqref{con0} and $K^\Sigma$, we may assume that $\V=\partial_z$.  As $\V^\perp$ is spacelike, the vector $w$ defining the rulings is necessarily spacelike. By a rotation in $\V^\perp$, we may further assume that $w=\partial_y$. Thus, the generating curve $\gamma$ is contained in the $xz$-plane, say $\gamma(s)=(x(s),0,z(s)),$  with $x'^2-z'^2=\varepsilon_1$. Consequently, $\Sigma$ can be parametrized as
\begin{equation}\label{gg}
\psi(s,t)=(x(s),t,z(s)),\qquad x'^2-z'^2=\varepsilon_1.
\end{equation}

\begin{corollary}\label{co44}
Let $\V=\partial_z$ and let $\Sigma$ be a cylindrical surface parametrized by \eqref{gg} with constant intrinsic sectional curvature $K^\Sigma=K_0$ with respect to the non-metric css-connection determined by $\V$. Then either $\Sigma$ is a plane or, up to translations of $s$ and $z$,
$$
z(s)=
\begin{cases}
\log\cosh(as)\ \text{ or }\ \log|\sinh(as)|,\quad a=\sqrt{2\varepsilon_1K_0}, & \varepsilon_1K_0>0,\\
\log|s|, & K_0=0,\\
\log|\cos(bs)|,\quad b=\sqrt{-2\varepsilon_1K_0}, & \varepsilon_1K_0<0,
\end{cases}
$$
where $x(s)= \pm\int^s\sqrt{\varepsilon_1+z'(u)^2}\,du$, defined on the intervals where $\varepsilon_1+z'^2>0$. If $\varepsilon_1=1$, this holds everywhere; otherwise ($\varepsilon_1=-1$) it requires $|z'|>1$. In particular, the branch
$z(s)=\log\cosh(as)$ admits timelike generating curves only if $a>1$, on intervals where $|a\tanh(as)|>1$. Similarly, for $\varepsilon_1=-1$, the branch $z(s)=\log|\sinh(as)|$ is defined everywhere if $a\geq1$, but only on intervals where $|a\coth(as)|>1$ if $a<1$, and the branch $z(s)=\log|\cos(bs)|$ only on intervals where $|b\tan(bs)|>1$.
\end{corollary}

\begin{proof}
Since $\V=\partial_z$, we have $p=\langle\V,\gamma'\rangle=-z'$. Integrating the expressions in \eqref{pRicc-solu} yields the stated formulas for $z$, up to suitable constants. The constant solutions $p=\pm a$ in \eqref{pRicc-solu}, and $p=0$ when $K_0=0$, give planes, see Theorem~\ref{clascyl}. Finally, $x(s)$ is recovered from $x'^2=\varepsilon_1+z'^2$.
\end{proof}

We next consider the cases in which $\V$ is spacelike or null, restricting initially to spacelike rulings. As before, we take $w=\partial_y$ and $\gamma(s)=(x(s),0,z(s))$.  When $\V$ is null, the choice $w=\partial_y$ entails no loss of generality, because the null rotations fixing $\V$ act transitively on the unit spacelike vectors in $\V^\perp$. When $\V$ is spacelike, the plane $\V^\perp$ is timelike and $w$ can also be timelike. This case is treated in Corollary~\ref{co410}.

\begin{corollary}\label{cor45}
Let $\V=\partial_x$ and let $\Sigma$ be a cylindrical surface parametrized by \eqref{gg} with constant intrinsic sectional curvature $K^\Sigma=K_0$ with respect to the non-metric css-connection determined by $\V$. Then, either $\Sigma$ is a plane or, up to translations of $s$ and $x$,
$$
x(s)=
\begin{cases}
-\log\cosh(as)\ \text{ or }\ -\log|\sinh(as)|,\quad a=\sqrt{2\varepsilon_1K_0}, & \varepsilon_1K_0>0,\\
-\log|s|, & K_0=0,\\
-\log|\cos(bs)|,\quad b=\sqrt{-2\varepsilon_1K_0}, & \varepsilon_1K_0<0,
\end{cases}
$$
where $z(s)=\pm\int^s\sqrt{x'(u)^2-\varepsilon_1}\,du$, defined on the intervals where $x'^2-\varepsilon_1>0$. If $\varepsilon_1=-1$, this holds everywhere; otherwise ($\varepsilon_1=1$) it requires $|x'|>1$. In particular, for $\varepsilon_1=1$, the branch $x(s)=-\log\cosh(as)$ is admissible only if $a>1$, on intervals where $|a\tanh(as)|>1$; the branch $x(s)=-\log|\sinh(as)|$ is defined everywhere if $a\geq1$, but only on intervals where $|a\coth(as)|>1$ if $a<1$; and the branch $x(s)=-\log|\cos(bs)|$ only on intervals where $|b\tan(bs)|>1$.
\end{corollary}

\begin{proof}
Since $\V=\partial_x$, we have $p=\langle\V,\gamma'\rangle=x'$. The constant solutions give planes. Integrating the expressions in \eqref{pRicc-solu} yields the stated formulas for $x$, up to suitable constants. Finally, $z(s)$ is recovered from $x'^2-z'^2=\varepsilon_1$.
\end{proof}

\begin{corollary}\label{co46}
Let $\V=\partial_x+\partial_z$ and let $\Sigma$ be a cylindrical surface parametrized by \eqref{gg} with constant intrinsic sectional curvature $K^\Sigma=K_0$ with respect to the non-metric css-connection determined by $\V$.  Then either $\Sigma$ is a plane or, up to a translation of $s$,
$$
x(s)=\frac12\int^s\left(p(u)+\frac{\varepsilon_1}{p(u)}\right)\,du,
\qquad
z(s)=\frac12\int^s\left(\frac{\varepsilon_1}{p(u)}-p(u)\right)\,du,
$$
defined on the intervals where $p$, given by one of the expressions in \eqref{pRicc-solu}, does not vanish.
\end{corollary}

\begin{proof}
Since $\V=\partial_x+\partial_z$, we have $p=\eta(\gamma')=x'-z'$. The arc-length condition $x'^2-z'^2=\varepsilon_1$ implies that $p\neq0$ everywhere. From the two identities it follows that $x'+z'=\varepsilon_1/p$, whence $x'=\tfrac12(p+\varepsilon_1/p)$ and $z'=\tfrac12(\varepsilon_1/p-p)$. Integrating the expressions for $x'$ and $z'$ gives the stated formulas.
\end{proof}

For non-metric connection, Corollaries~\ref{co44}, \ref{cor45}, \ref{co46} and \ref{co410} make explicit  Theorem~\ref{clascyl} for cylindrical surfaces whose rulings are orthogonal to $\V$, according to the causal character of $\V$ and of the rulings. Corollary~\ref{co44} is the Lorentzian, non-metric analogue of Corollary 3.3 in \cite{alm1}, which gives the corresponding classification in the Euclidean case, with the same trichotomy according to the sign of $K$.

\subsection{The metric case ($\alpha=1$)}

We now consider the metric css-connection, corresponding to $\alpha=1$. Since the rulings are assumed to be orthogonal to $\V$, we have $c=\langle\V,w\rangle=0$. Thus, by Theorem \ref{clascyl}(iii), the condition $K^\Sigma=K_0$ implies
\begin{equation}\label{linear-solu}
p(s)=-\varepsilon_1K_0s.
\end{equation}
As in the previous subsection, we investigate the different causal characters of $\V$ separately, taking $w=\partial_y$ and $\gamma(s)=(x(s),0,z(s))$.

\begin{corollary}\label{co47}
Let $\V=\partial_z$ and let $\Sigma$ be a cylindrical surface parametrized by \eqref{gg} with constant intrinsic sectional curvature $K^\Sigma=K_0$ with respect to the metric css-connection determined by $\V$. If $K_0=0$, then $\Sigma$ is a plane. If $K_0\neq0$, then, up to translations of $s$ and $z$, $z(s)=\tfrac12\varepsilon_1K_0s^2$, and
$$
x(s)=\pm
\begin{cases}
\dfrac{s}{2}\sqrt{1+K_0^2s^2}+\dfrac{1}{2|K_0|}\sinh^{-1}(|K_0|s), & \varepsilon_1=1,\\
\dfrac{s}{2}\sqrt{K_0^2s^2-1}-\dfrac{1}{2|K_0|}\cosh^{-1}(|K_0|s), & \varepsilon_1=-1,\ s\geq 1/|K_0|,
\end{cases}
$$
the case $s\leq-1/|K_0|$ being symmetric.
\end{corollary}

\begin{proof}
Since $\V=\partial_z$, we have $p=\langle\V,\gamma'\rangle=-z'$. The constant solutions for $z'$ give planes. Otherwise, integrating \eqref{linear-solu} yields, up to translations of $z$, $z(s)=\tfrac12\varepsilon_1K_0s^2$. Finally, $x(s)$ is recovered from $x'^2=\varepsilon_1+z'^2=\varepsilon_1+K_0^2s^2$, yielding the two cases above.
\end{proof}

The following result is obtained by the same argument as     in the proof of Corollary~\ref{co47}.

\begin{corollary}\label{co48}
Let $\V=\partial_x$ and let $\Sigma$ be a cylindrical surface parametrized by \eqref{gg} with constant intrinsic sectional curvature $K^\Sigma=K_0$ with respect to the metric css-connection determined by $\V$. If $K_0=0$, then $\Sigma$ is a plane. If $K_0\neq0$, then, up to translations of $s$ and $x$, $x(s)=-\tfrac12\varepsilon_1K_0s^2$, and
$$
z(s)=  \pm 
\begin{cases}
\dfrac{s}{2}\sqrt{K_0^2s^2-1}-\dfrac{1}{2|K_0|}\cosh^{-1}(|K_0|s), & \varepsilon_1=1,\ s\geq1/|K_0|,\\
\dfrac{s}{2}\sqrt{1+K_0^2s^2}+\dfrac{1}{2|K_0|}\sinh^{-1}(|K_0|s), & \varepsilon_1=-1,
\end{cases}
$$
the case $s\leq-1/|K_0|$ being symmetric.
\end{corollary}

\begin{corollary}\label{co49}
Let $\V=\partial_x+\partial_z$ and let $\Sigma$ be a cylindrical surface parametrized by \eqref{gg} with constant intrinsic sectional curvature $K^\Sigma=K_0$ with respect to the metric css-connection determined by $\V$. If $K_0=0$, then $\Sigma$ is a plane. If $K_0\neq0$, then, up to translations of $s,x,z$, 
$$
x(s)=-\frac{\varepsilon_1K_0}{4}s^2-\frac{1}{2K_0}\log|s|,
\qquad
z(s)=\frac{\varepsilon_1K_0}{4}s^2-\frac{1}{2K_0}\log|s|,\qquad s\neq0.
$$
\end{corollary}

\begin{proof}
If $K_0=0$, then $p=x'-z'$ is a nonzero constant. Hence, both $x$ and $z$ are linear functions, so $\Sigma$ is a plane. Assume that $K_0\neq0$. From $x'-z'=p$ and $x'+z'=\varepsilon_1/p=-1/(K_0s)$, we obtain $x'=\tfrac12\big({-\varepsilon_1K_0s}-\tfrac{1}{K_0s}\big)$ and $z'=\tfrac12\big(\varepsilon_1K_0s-\tfrac1{K_0s}\big)$, which integrate to the stated formulas.
\end{proof}

\subsection{Timelike  and null rulings}

Suppose now that  $\V$ is spacelike. Then the plane $\V^\perp$ is timelike and the rulings orthogonal to $\V$ may be timelike. Without loss of generality, assume that $\V=\partial_x$ and $w=\partial_z$, with $\varepsilon_2=-1$. Hence, $\gamma$ is contained in the spacelike $xy$-plane, say $\gamma(s)=(x(s),y(s),0)$, with $x'^2+y'^2=1$. Consequently, $\Sigma$ is timelike and can be parametrized as
\begin{equation}\label{ggt}
\psi(s,t)=(x(s),y(s),t),\qquad x'^2+y'^2=1.
\end{equation}

\begin{corollary}\label{co410}
Let $\V=\partial_x$ and let $\Sigma$ be a non-planar cylindrical surface parametrized by \eqref{ggt} with constant intrinsic sectional curvature $K^\Sigma=K_0$ with respect to the css-connection $\nabla$ determined by $\V$.
\begin{enumerate}
\item[(i)] If $\nabla$ is non-metric, then, up to translations of $s$ and $x$,
$$
x(s)=
\begin{cases}
-\log\cosh(as)\ \text{ or }\ -\log|\sinh(as)|,\quad a=\sqrt{2K_0}, & K_0>0,\\
-\log|s|, & K_0=0,\\
-\log|\cos(bs)|,\quad b=\sqrt{-2K_0}, & K_0<0,
\end{cases}
$$
and $y(s)=\pm\int^s\sqrt{1-x'(u)^2}\,du$, defined on the intervals where $x'^2<1$. The branch $-\log|\sinh(as)|$ requires $a<1$. 
\item[(ii)] If $\nabla$ is metric, then $K_0\neq0$ and, up to translations of $s$, $x$ and $y$,
$$
x(s)=-\frac{K_0}{2}s^2,\qquad y(s)=\pm\Big(\frac{s}{2}\sqrt{1-K_0^2s^2}+\frac{1}{2K_0}\sin^{-1}(K_0s)\Big),\qquad |K_0s|<1.
$$
\end{enumerate}
\end{corollary}

\begin{proof}
Since $\V=\partial_x$, we have $p=x'$. In the non-metric case, integrating the non-constant solutions in \eqref{pRicc-solu} gives the stated expressions for $x$. In the metric case, the non-planarity assumption implies $K_0\neq0$ and integrating \eqref{linear-solu} yields $x(s)=-K_0s^2/2$. Finally, the arc-length condition
$x'^2+y'^2=1$ gives $y'=\pm\sqrt{1-x'^2},$ from which the stated expressions for $y$ follow.
\end{proof}

Finally, we consider cylindrical surfaces with null rulings. These surfaces are not covered by \eqref{eqcyl}, since a null ruling vector cannot be normalized and its orthogonal complement is degenerate.

\begin{proposition}\label{nullrulings}
Let $\Sigma$ be a nondegenerate cylindrical surface parametrized by $\psi(s,t)=\gamma(s)+t\,w$, where $w$ is a fixed null vector and $\gamma$ is a curve transversal to $w$. Then $\Sigma$ is timelike, $G=H^0=0$ and 
$$
K^\Sigma=\frac{(1-\alpha)^2}{2}\big(\delta-\eta(N)^2\big),
$$
where $K^\Sigma$ is the intrinsic sectional curvature of $\Sigma$ with respect to the css-connection determined by $\V$. Consequently, if $\Sigma$ is non-planar, the following statements hold:
\begin{enumerate}
\item[(i)] For the metric css-connection, $K^\Sigma=0$.
\item[(ii)]  For the non-metric css-connection, $K^\Sigma$ is constant if and only if $\langle\V,w\rangle=0$, and in this case $K^\Sigma=0$.
\end{enumerate}
\end{proposition}

\begin{proof}
Since $\Sigma$ is nondegenerate and contains the null direction $w$, it is timelike, so $\varepsilon=1$. As $\psi_t=w$ and $\psi_{st}=\psi_{tt}=0$, we have $m^0=n^0=0$ and so $G=0$. Moreover, $g=\langle w,w\rangle=0$; hence $H^0=\frac{\varepsilon}{2}\frac{l^0g-2m^0f+n^0e}{eg-f^2}=0$. Proposition~\ref{pr31} gives the expression of $K^\Sigma$, and case (i) follows.

For case (ii), $N$ is orthogonal to $w$, so $N$ lies in the degenerate plane $w^\perp$, which is spanned by $w$ and a unit spacelike vector $e$. Since $\langle N,N\rangle=1$, we have $N=\pm e+\lambda w$ for some function $\lambda$, with constant sign by continuity, and $\eta(N)=\pm\langle\V,e\rangle+\lambda\langle\V,w\rangle$. Now $K^\Sigma$ is constant if and only if $\eta(N)$ is constant. If $\langle\V,w\rangle\neq0$, then $\eta(N)$ is constant if and only if $\lambda$ is constant, that is, if and only if  $N$ is constant, in which case $\Sigma$ is a plane. If $\langle\V,w\rangle=0$, then $\V\in w^\perp$, so $\V=ae+bw$ with $\delta=a^2$ and $\eta(N)=\pm a$. Hence $K^\Sigma=\frac12(a^2-a^2)=0$.
\end{proof}

The timelike translation surfaces $z=f(x)\pm y$ considered in Example~\ref{translation} of Section~\ref{sec3}
are cylindrical surfaces with null rulings parallel to $(0,1,\pm1)$. Since $\langle\partial_z,(0,1,\pm1)\rangle\neq0$,
Proposition~\ref{nullrulings} implies that their intrinsic sectional curvature $K^\Sigma=K$ with respect to the non-metric
css-connection cannot be constant unless the surfaces are planes.

Corollary~\ref{cor32}(i), Theorem~\ref{clascyl} and Proposition~\ref{nullrulings} together classify all cylindrical surfaces of $\R_1^3$ with constant intrinsic sectional curvature with respect to css-connections, covering all causal characters of $\V$ and of the ruling direction.

\section{Rotational surfaces with constant sectional curvature with respect to css-connections} \label{sec5}

In $\R_1^3$, one-parameter groups of linear isometries that fix pointwise a straight line through the origin (called the axis) are of three types according to the causal character of the axis: elliptic rotations for a timelike axis, hyperbolic rotations for a spacelike axis and parabolic rotations for a null axis. A rotational surface is a surface invariant under one of these groups.

If $\Phi$ is a linear isometry of $\R_1^3$, then $\Phi$ maps the css-connection determined by $\V$ onto the css-connection determined by $\Phi(\V)$. Hence a group of rotations preserves the css-connection if and only if it fixes $\V$, that is, if and only if $\V$ is parallel to the axis. In this case $K^\Sigma$ is constant along the orbits and the condition $K^\Sigma=K_0$ becomes an ordinary differential equation for the generating curve. Throughout this section we assume that $\V$ is parallel to the axis. The other case is discussed in Remark~\ref{rem-nonaxis}. Up to a linear isometry, and by Remark~\ref{rem-norm}, there are four cases. In each of them,   the profile curve $\gamma$ is parametrized by arc length, $\langle\gamma',\gamma'\rangle=\varepsilon_1$, and $\theta\in\R$:
\begin{enumerate}
\item[(E)] Timelike axis, $\V=\partial_z$: $\Psi(s,\theta)=(x(s)\cos\theta,x(s)\sin\theta,z(s))$, $x>0$.
\item[(H1)] Spacelike axis, $\V=\partial_x$: $\Psi(s,\theta)=(x(s),z(s)\sinh\theta,z(s)\cosh\theta)$, $z>0$.
\item[(H2)] Spacelike axis, $\V=\partial_x$: $\Psi(s,\theta)=(x(s),y(s)\cosh\theta,y(s)\sinh\theta)$, $y>0$.
\item[(P)] Null axis, $\V=(0,1,1)$: $\Psi(s,\theta)=A(\theta)\gamma(s)$, $\gamma(s)=(0,y(s),z(s))$, $y-z>0$, where
$$
A(\theta)=\begin{pmatrix}1&-\theta&\theta\\ \theta&1-\frac{\theta^2}{2}&\frac{\theta^2}{2}\\ \theta&-\frac{\theta^2}{2}&1+\frac{\theta^2}{2}\end{pmatrix}.
$$
\end{enumerate}
In the hyperbolic case (spacelie axis), the orbits of the points with $|z|>|y|$ meet the plane $y=0$, which gives (H1). The orbits of the points with $|y|>|z|$ meet the plane $z=0$, which gives (H2). The points of the null planes $y=\pm z$ have null orbits and are not considered. In case (P), the matrices $A(\theta)$ are the parabolic rotations fixing the null vector $(0,1,1)$, and the orbit of every point with $y\neq z$ meets the plane $x=0$ exactly once. By Remark~\ref{rem-norm}, any two nonzero null vectors are related by a Lorentz isometry, so the choice $\V=(0,1,1)$ in case (P), instead of $\V=\partial_x+\partial_z$ as in Section~\ref{sec4}, is only a matter of convenience.

In each case we define two functions $\rho(s)$ and $\xi(s)$, and a sign $\sigma$, as follows:
$$
\begin{array}{c|c|c|c|c}
 & \rho & \xi & \sigma & \varepsilon_1\\ \hline
\text{(E)} & x & -z & 1 & x'^2-z'^2\\
\text{(H1)} & z & x & -1 & x'^2-z'^2\\
\text{(H2)} & y & x & 1 & x'^2+y'^2=1\\
\text{(P)} & y-z & y-z & 1 & y'^2-z'^2
\end{array}
$$
In all cases the first fundamental form of $\Sigma$ is
$$
\varepsilon_1\,ds^2+\sigma\rho(s)^2\,d\theta^2,
$$
so $\Sigma$ is spacelike if $\varepsilon_1=\sigma=1$ and timelike otherwise. Moreover, since $\V$ is fixed by the rotations, the height function is $\langle\V,\Psi(s,\theta)\rangle=\langle A(\theta)^{-1}\V,\gamma(s)\rangle=\langle\V,\gamma(s)\rangle=\xi(s)$. Here $A(\theta)$ denotes the rotation of the corresponding type.

\begin{proposition}\label{throt}
Let $\Sigma$ be a rotational surface as above. Then its intrinsic sectional curvature with respect to the css-connection is
\begin{equation}\label{Krot}
K^\Sigma=\varepsilon_1\left(\frac{(1-\alpha)^2}{2}\,\xi'^2-\frac{\rho''}{\rho}-\frac{1+\alpha}{2\rho}\,(\rho \xi')'\right).
\end{equation}
\end{proposition}

\begin{proof}
For the metric $\varepsilon_1ds^2+\sigma\rho^2d\theta^2$ and a function $\xi=\xi(s)$, we compute
$$
G=-\varepsilon_1\frac{\rho''}{\rho},\qquad \langle\nabla^{0_\Sigma}\xi,\nabla^{0_\Sigma}\xi\rangle=\varepsilon_1\xi'^2,\qquad \Delta^\Sigma \xi=\varepsilon_1\frac{(\rho \xi')'}{\rho}.
$$
The result follows from Corollary~\ref{intrinsic}.
\end{proof}

\begin{corollary}\label{rotexamples}
Let $\Sigma$ be a rotational surface as above. Then:
\begin{enumerate}
\item[(i)] In cases (E), (H1) and (H2), the planes orthogonal to the axis ($\xi'=0$) have $K^\Sigma=0$.
\item[(ii)] In cases (E), (H1) and (H2), the surfaces with $\rho$ constant are cylindrical surfaces with rulings parallel to $\V$, and $K^\Sigma=\delta\frac{(1-\alpha)^2}{2}$.
\item[(iii)] In case (P), neither  of the these occurs.
\end{enumerate}
\end{corollary}

\begin{proof}
(i) If $\xi'=0$, the arc-length condition gives $\rho=\pm s$ up to a translation, so $\rho''=0$ and $K^\Sigma=0$ by \eqref{Krot}. (ii) If $\rho'=0$, then $\xi'^2=1$ and $\varepsilon_1=\delta$, so \eqref{Krot} gives $K^\Sigma=\varepsilon_1\frac{(1-\alpha)^2}{2}=\delta\frac{(1-\alpha)^2}{2}$, in agreement with Corollary~\ref{cor32}(ii). (iii) In case (P), $\rho=\xi$, and $\rho'=0$ would give $\varepsilon_1=y'^2-z'^2=(y'-z')(y'+z')=0$.
\end{proof}

\subsection{The metric case ($\alpha=1$)}

For $\alpha=1$, \eqref{Krot} becomes $\varepsilon_1K^\Sigma\rho=-(\rho'+\rho \xi')'$. Hence:

\begin{corollary}\label{rotmetric}
A rotational surface as above has constant intrinsic sectional curvature $K_0$ with respect to the metric css-connection if and only if
$$
(\rho'+\rho \xi')'=-\varepsilon_1K_0\,\rho.
$$
In particular, $K_0=0$ if and only if $\rho'+\rho \xi'=c_1$ for some constant $c_1$.
\end{corollary}

In the case $K_0=0$, the equation $\rho'+\rho \xi'=c_1$ can be integrated. On an interval where $\rho'\neq0$, we write the profile curve as a graph over $r=\rho$: $z=u(x)$ in case (E), $x=u(z)$ in case (H1) and $x=u(y)$ in case (H2). Let $q=u'(r)$.

\begin{corollary}\label{rotflat}
A rotational surface of type (E), (H1) or (H2) with $\rho'\neq0$ has vanishing intrinsic sectional curvature with respect to the metric css-connection if and only if there is a constant $\lambda$ such that
\begin{enumerate}
\item[(E)] $(rq-1)^2=\lambda(1-q^2)$, that is, $q=\dfrac{r\pm\sqrt{\lambda(r^2+\lambda-1)}}{r^2+\lambda}$.
\item[(H1)] $(rq+1)^2=\lambda(q^2-1)$, that is, $q=\dfrac{-r\pm\sqrt{\lambda(\lambda+1-r^2)}}{r^2-\lambda}$.
\item[(H2)] $(rq+1)^2=\lambda(1+q^2)$, $\lambda\geq0$, that is, $q=\dfrac{-r\pm\sqrt{\lambda(r^2+1-\lambda)}}{r^2-\lambda}$.
\end{enumerate}
The profile curve is then obtained from $u(r)=\int^r q(\tau)\,d\tau$. In particular, for $\lambda=0$ we obtain the surfaces $z=\log x$ in case (E), $x=-\log z$ in case (H1) and $x=-\log y$ in case (H2). If $\rho'=0$, then $\Sigma$ is one of the cylindrical surfaces of Corollary~\ref{rotexamples}(ii), which also have $K^\Sigma=0$ for $\alpha=1$.
\end{corollary}

\begin{proof}
In case (E), we have $\rho=r$, $\xi=-u(r)$ and $\varepsilon_1ds^2=(1-q^2)dr^2$. Hence, $\rho'=dr/ds=\pm|1-q^2|^{-1/2}$ and $\xi'=-q\rho'$. Therefore, $\rho'+\rho \xi'=\rho'(1-rq)=c_1$ gives $(rq-1)^2=c_1^2|1-q^2|$, which is (E) with $\lambda=c_1^2\operatorname{sign}(1-q^2)$. Cases (H1) and (H2) are analogous, now with $\xi'=q\rho'$ and, respectively, $\varepsilon_1ds^2=(q^2-1)dr^2$ and $ds^2=(1+q^2)dr^2$. These give $(1+rq)^2=c_1^2|q^2-1|$ and $(1+rq)^2=c_1^2(1+q^2)$, respectively. Solving the quadratic equations in $q$ gives the explicit expressions.
\end{proof}

\subsection{The null axis}

When the axis is null, \eqref{Krot} can be integrated explicitly in the metric and non-metric cases and for every $K_0$. Since $\varepsilon_1=(y'-z')(y'+z')$, the function $\omega=y-z$ satisfies $\omega'\neq0$, so we may write the profile curve as $y+z=u(\omega)$, $\omega>0$. Then
$$
\Psi(\omega,\theta)=\Big(-\theta \omega,\ \frac{u+\omega}{2}-\frac{\theta^2}{2}\omega,\ \frac{u-\omega}{2}-\frac{\theta^2}{2}\omega\Big),
$$
where the first fundamental form is $u'(\omega)\,d\omega^2+\omega^2\,d\theta^2$, and $\Sigma$ is spacelike (resp. timelike) if $u'>0$ (resp. $u'<0$).

\begin{theorem}\label{thnull}
Let $\Sigma$ be a rotational surface with null axis parametrized as above. Then its intrinsic sectional curvature with respect to the css-connection is given by
$$
K^\Sigma=\frac{2(\omega-1)u'+(\omega+2)u''}{4\omega\,u'^2}\quad(\alpha=0),\qquad K^\Sigma=\frac{(\omega+1)u''-2u'}{2\omega\,u'^2}\quad(\alpha=1).
$$
Therefore $K^\Sigma$ is constant, $K^\Sigma=K_0$, if and only if there is a constant $C$ such that
\begin{enumerate}
\item[(i)] for the non-metric connection,
$$
u'(\omega)=\frac{2(\omega+2)^6}{C\,e^{2\omega}+K_0\,P(\omega)},
$$
where $P(\omega)=4\omega^6+52\omega^5+290\omega^4+900\omega^3+1670\omega^2+1798\omega+899$.
\item[(ii)] for the metric connection,
$$
u'(\omega)=\frac{3(\omega+1)^2}{C-K_0\,\omega^2(2\omega+3)}.
$$
\end{enumerate}
The profile curve is obtained from $u(\omega)=\int^\omega u'(\tau)\,d\tau$, on the intervals where the denominators do not vanish. In particular, the surfaces with $K_0=0$ are given by $u'=\frac{2}{C}(\omega+2)^6e^{-2\omega}$ for the non-metric connection and by $u=\frac{(\omega+1)^3}{C}$ for the metric connection, up to an additive constant.
\end{theorem}

\begin{proof}
The expression for $K^\Sigma$ follows from Corollary~\ref{intrinsic} applied to the metric $u'd\omega^2+\omega^2d\theta^2$ and the height function $\xi=\omega$, or, equivalently, from \eqref{Krot}. Set $v=1/u'$, so that $v'=-u''/u'^2$.

For $\alpha=0$, $K^\Sigma=K_0$ is equivalent to $(\omega+2)u''=4K_0\omega u'^2-2(\omega-1)u'$, that is,
$$
v'-\frac{2(\omega-1)}{\omega+2}\,v=-\frac{4K_0\omega}{\omega+2}.
$$
An integrating factor is $e^{-2\omega}(\omega+2)^6$, and $\int \omega(\omega+2)^5e^{-2\omega}\,d\omega=-\frac18P(\omega)e^{-2\omega}$. Hence $e^{-2\omega}(\omega+2)^6v=c_1+\frac{K_0}{2}P(\omega)e^{-2\omega}$, which gives (i) with $C=2c_1$.

For $\alpha=1$, $K^\Sigma=K_0$ is equivalent to $(\omega+1)u''=2u'+2K_0\omega u'^2$, that is, $v'+\frac{2}{\omega+1}v=-\frac{2K_0\omega}{\omega+1}$. Then $\big((\omega+1)^2v\big)'=-2K_0\omega(\omega+1)$, so $(\omega+1)^2v=c_1-K_0\big(\frac23\omega^3+\omega^2\big)$, which gives (ii) with $C=3c_1$.
\end{proof}

\subsection{The remaining cases}

For the axes (E), (H1) and (H2) with $\alpha=0$, or with $\alpha=1$ and $K_0\neq0$, if we write the profile curve as a graph, as in Corollary~\ref{rotflat}, the condition $K^\Sigma=K_0$ is equivalent to the following first-order differential equation for $q=u'(r)$, where $E$ denotes the coefficient of $dr^2$ in the first fundamental form:
\begin{itemize}
\item[(E)] $E=1-q^2$: \ $2K_0rE^2=Eq(rq+1)+(r-2q)q'$ if $\alpha=0$, and $K_0rE^2=Eq+(r-q)q'$ if $\alpha=1$.
\item[(H1)] $E=q^2-1$: \ $2K_0rE^2=Eq(rq-1)+(r+2q)q'$ if $\alpha=0$, and $K_0rE^2=-Eq+(r+q)q'$ if $\alpha=1$.
\item[(H2)] $E=1+q^2$: \ $2K_0rE^2=Eq(rq-1)+(2q-r)q'$ if $\alpha=0$, and $K_0rE^2=-Eq+(q-r)q'$ if $\alpha=1$.
\end{itemize}
By the existence and uniqueness theorem for ordinary differential equations, for every $K_0$ and every initial condition $q(r_0)=q_0$ with $E(q_0)\neq0$ and nonzero coefficient of $q'$, there is a unique local solution. Hence, for every $K_0\in\R$ there is a two-parameter family of rotational surfaces of each type with constant sectional curvature $K_0$, determined by the initial value $q_0$ and a translation along the axis. Apart from the solutions of Corollaries~\ref{rotexamples} and \ref{rotflat}, we do not know of any closed-form solutions of these equations.

\begin{remark}\label{rem-nonaxis}
If $\V$ is not parallel to the axis, then the rotations do not preserve the css-connection and the height function $\xi$ depends on $\theta$, so that the condition $K^\Sigma=K_0$ becomes a partial differential equation in $(s,\theta)$. This case is not considered here.
\end{remark}

\section*{Acknowledgements}
Rafael L\'opez has been partially supported by MINECO/MICINN/FEDER grant no.
PID2023-150727NB-I00, and by the ``Mar\'{\i}a de Maeztu'' Excellence Unit IMAG, reference CEX2020-001105-M, funded by MICINN/AEI/10.13039/501100011033.



\begin{thebibliography}{99}

\bibitem{ag} N. S. Agashe, M. R. Chafle, A semi-symmetric non-metric connection on a Riemannian manifold. Indian J. Pure Appl. Math. 23 (1992), 399--409.

\bibitem{ac} N. S. Agashe, M. R. Chafle, On submanifolds of a Riemannian manifold with a semi-symmetric non-metric connection. Tensor 55 (1994), 120--130.

\bibitem{alm0} M.E. Ayd{\i}n, R. L\'{o}pez, A. Mihai, Classification of translation surfaces in $\mathbb{R}^3$ with constant sectional curvature for the canonical snm-connection. Bull. Malays. Math. Sci. Soc. 48 (2025), 114. 

\bibitem{alm1} M.E. Ayd{\i}n, R. L\'{o}pez, A. Mihai, Constant sectional curvature surfaces with a semi-symmetric non-metric connection. J. Math. Anal. Appl. 542 (2025), 128795. 


\bibitem{ds} K. Duggal, R. Sharma, Semi-symmetric metric connections in a semi-Riemannian manifold. Indian J. Pure Appl. Math. 17 (1986), 1276--1283.

\bibitem{fs} A. Friedmann, J. A. Schouten, \"{U}ber die Geometrie der halbsymmetrischen \"{U}bertragungen. Math. Z. 21 (1924), 211--223. 


\bibitem{hay} H. Hayden, Subspaces of a space with torsion, Proc. London Math. Soc. s2-34 (1932), 27--50. 

\bibitem{im} T. Imai, Hypersurfaces of a Riemannian manifold with semi-symmetric metric connection. Tensor (N. S.) 23 (1972), 300--306. 

\bibitem{lyl} C. W. Lee, D. W. Yoon, J. W. Lee, Optimal inequalities for the Casorati curvatures of submanifolds of real space forms endowed with semi-symmetric metric connections. J. Inequal. Appl. 2014 (2014), 327. 

\bibitem{mm0} A. Mihai, I. Mihai, A note on a well-defined sectional curvature of a semi-symmetric non-metric connection. Int. Electron. J. Geom. 17 (2024), 15--23. 

 \bibitem{mc1} A. Mihai, C. \"{O}zg\"{u}r, Chen inequalities for submanifolds of real space forms with a semi-symmetric metric connection. Taiwanese J. Math. 14 (2010), 1465--1477. 



\bibitem{mo} I. Mihai, A. Olteanu, Chen-like inequalities for submanifolds in Kähler manifolds admitting semi-symmetric non-metric connections. Symmetry, 16 (2024), 1401. 

\bibitem{mmme} M. Mohammed, F. Massamba, I. Mihai, A.E.A.M.A. Elamin, M.S. Aldien, Some Chen inequalities for submanifolds in trans-Sasakian manifolds admitting a semi-symmetric non-metric connection. Axioms 13 (2024), 195. 

\bibitem{mma} M. Mohammed, I. Mihai, A. Olteanu, Pinching results for submanifolds in Lorentzian–Sasakian manifolds endowed with a semi-symmetric non-metric connection. Mathematics 12 (2024), 3651. 

\bibitem{na} Z. Nakao, Submanifolds of a Riemannian manifold with semisymmetric metric connections. Proc. Amer. Math. Soc. 54 (1976), 261--266. 


 
\bibitem{tr} M. M. Tripathi, A new connection in a Riemannian manifold. Int. Electron. J. Geom. 1 (2008), 15--24. 

\bibitem{wa} Y. Wang,  Minimal translation surfaces with respect to semi-symmetric connections in $\mathbb{R}^3$ and $\mathbb{R}^3_1$. Bull. Korean Math. Soc. 58 (2021),  959--972. 



\bibitem{ya} K. Yano, On semi symmetric metric connection. Rev. Roum. Math. Pures Appl. 15 (1970), 1579--1591.

\bibitem{yay0} A. Y\"{u}cesan, N. Ayy{\i}ld{\i}z, Non-degenerate hypersurfaces of a semi–Riemannian manifold with a semi–symmetric non-metric connection. Math. Rep. (Bucur) 14 (2012), 209--219.

\bibitem{yay1} A. Y\"{u}cesan, N. Ayy{\i}ld{\i}z, Non-degenerate hypersurfaces of a semi-Riemannian manifold with a semi-symmetric metric connection. Arch. Math. (Brno) 44 (2008), 77--88.




\bibitem{zz} P. Zhang, L. Zhang, Casorati inequalities for submanifolds in a Riemannian manifold of quasi-constant curvature with a semi-symmetric metric connection. Symmetry 8(4) (2016), 19. 

\end{thebibliography}
 \end{document}